\documentclass[11pt,a4paper]{article}
 \usepackage{euscript}
\usepackage{amsmath}
\usepackage{graphicx}
\usepackage{amsthm}
\usepackage{amssymb}
\usepackage{epsfig}
\usepackage{url}
\usepackage{colonequals}
\usepackage{amsmath,amscd}
\theoremstyle{theorem}
\newtheorem{theorem}{Theorem}[section]
\newtheorem{corollary}[theorem]{Corollary}

\newtheorem{lemma}[theorem]{Lemma}
\newtheorem{proposition}[theorem]{Proposition}

\newtheorem{definition}{Definition}[section]

\newtheorem{remark}{Remark}[section]
\numberwithin{equation}{section}

\usepackage{tikz}
\usetikzlibrary{arrows.meta, positioning}

\DeclareSymbolFont{symbolsC}{U}{txsyc}{m}{n}
\DeclareMathSymbol{\notniFromTxfonts}{\mathrel}{symbolsC}{61}
\title{The Central Nullstellensatz over Centrally Algebraically Closed Division Rings} 
\author{Masood Aryapoor\\
	\tiny{\textit{Division of Mathematics and Physics}}\\
	\tiny{\textit{M\"{a}lardalen  University}}\\
	\tiny{\textit{Hamngatan 15, 632 17, Eskilstuna, 
			Sweden
	}}
}
 \date{}
\begin{document}
 \maketitle

\begin{abstract}

We introduce the concept of centrally algebraically closed division rings and show that a division ring satisfies the central Nullstellensatz if and only if it is centrally algebraically closed. We also show that every division ring can be embedded in a centrally algebraically closed division ring. 

\end{abstract}
\begin{section}{Introduction} 

In commutative algebra, the weak Nullstellensatz states that for an algebraically closed field $k$,  every maximal ideal of the ring $k[x_1,\dots,x_n]$ is of the form $(x_1-a_1,\dots,x_n-a_n)$ for some point $(a_1,\dots,a_n)\in k^n$. One can prove the weak Nullstellensatz using Zariski's lemma, which states that every ideal of $k[x_1,\dots,x_n]$, where $k$ is a (not necessarily algebraically closed) field, has finite codimension   as a vector space over $k$. In \cite{AmitsurSmall}, Amitsur and Small proved a noncommutative version of Zariski's lemma: 
\begin{theorem}[Amitsur--Small Theorem]\label{thm: Weak Nullstellensatz over division rings}\label{thm: amitsur small}
	Let $D$ be a division ring, and let $D[x_1,\dots,x_n]$ denote the polynomial ring  in $n$ central indeterminates over $D$. Then every simple (left or right) $D[x_1,\dots,x_n]$-module is finite-dimensional as a vector space over $D$.
\end{theorem}            
The principal aim of this work is to leverage the Amitsur--Small Theorem to establish an analogue of the weak Nullstellensatz for ``algebraically closed division rings". Before stating our main result, it is worth noting that the literature contains several logically nonequivalent definitions of “algebraically closed division rings”.  The notion most pertinent to our study is that of a \textit{right algebraically closed division} ring: A division ring $D$ is called \textit{right algebraically closed} if every equation $x^n+a_{n-1}x^{n-1}+\cdots+a_1x+a_0=0$, where $a_i\in D$, has a solution in $D$.   A classical result of Niven and Jacobson states that the division ring of quaternions over the field of real numbers is right algebraically closed. In \cite{AlonNullestellensatz}, Alon and Paran proved the following weak form of the Nullstellensatz over $\mathbb{H}$, henceforth referred to as the \textit{central Nullstellensatz}:    
\begin{theorem}[Alon--Paran Theorem]\label{thm: Weak Nullstellensatz over quaternions}
	Let $\mathbb{H}$ be the division ring of quaternions over $\mathbb{R}$, and let $\mathbb{H}[x_1,\dots,x_n]$ denote the polynomial ring  in $n$ central indeterminates over $\mathbb{H}$. Then a left ideal $I$ of   $\mathbb{H}[x_1,\dots,x_n]$ is maximal if and only if it is generated by $x_1-a_1,\dots,x_n-a_n$ for some point $(a_1,\dots,a_n)\in \mathbb{H}^n$ with pairwise commuting components.  
\end{theorem} 
The question of whether the central Nullstellensatz holds over every right algebraically closed division ring was raised by the author in \cite{aryapoor2024explicit}. The same question was also raised in \cite{chapman2025amitsur}, where Chapman and Paran called a division ring satisfying the central Nullstellensatz a \textit{Nullstellensatz ring} \cite[Definition 4.1]{chapman2025amitsur}. Our main result is that a division ring $D$ is a Nullstellensatz ring if and only if for every $(a_1,\dots,a_n)\in D^n$ with pairwise commuting components, the centralizer of $\{a_1,\dots,a_n\}$ in $D$ is right algebraically closed.  Furthermore, we present a version of the strong Nullstellensatz for such division rings, which is in line with \cite{aryapoor2024explicit}. We also study centralizers and their behavior regarding left algebraicity.

The paper is organized as follows. Section \ref{sec: prelimiaries} provides the necessary preliminaries. In Section \ref{sec: algebraic elements}, we discuss the notion of left algebraic elements and related notions in the context of division rings.  Section \ref{sec: null} contains the main result of the paper. The final section is devoted to some further results and open problems.


\end{section}
\begin{section}{Preliminaries} \label{sec: prelimiaries}
This section presents essential preliminaries, most of which are established results. For more details, we refer the reader to Lam's book \cite[Section 5.16]{lam2001first}.


	\begin{subsection}{Evaluation of polynomials over division rings}\label{subsec: evaluation}
		Let $D$ be a division ring, and let $D[x]$ be the ring of polynomials over $D$ in a central indeterminate $x$. Any nonzero element of $D[x]$ can be uniquely expressed as as a (left) polynomial  $a_nx^n+\cdots+a_1x+a_0$, where $a_0,\dots,a_n\in D$ with $a_n\neq 0$. The (left) \textit{evaluation} of  a polynomial  $p(x)=a_nx^n+\cdots+a_1x+a_0\in D[x]$ at $a\in D$ is defined to be
		\[
		p(a)\colonequals a_na^n+\cdots+a_1a+a_0. 
		\]
		The \textit{right evaluation} of  a polynomial  $p(x)=a_nx^n+\cdots+a_1x+a_0\in D[x]$ at $a\in D$ is defined to be
		\[
		p_r(a)\colonequals a^na_n+\cdots+aa_1+a_0. 
		\]
		We have the product formula 
		\begin{equation}\label{equ: left product formula in one variable}
			(pq)(a)=
			\begin{cases}
				0& \text{ if } q(a)=0,\\
				p\left( q(a)\,a\,q(a)^{-1} \right)q(a) & \text{ if  } q(a)\neq 0,
			\end{cases}
		\end{equation}
		where $p,q\in D[x]$ and $a\in D$.  
		An element $a\in D$ is called a \textit{right root} (respectively, left root)   of a polynomial $p(x)\in D[x]$ if $p(a)=0$ (respectively, $p_r(a)=0$). It is easy to see that $p(a)=0$ if and only if $p(x)\in D[x](x-a)$.  
		
		Next, we recall some well-known facts about roots of left polynomials over division rings (see \cite[Section 3]{lam2008wedderburn}). Let $p(x)\in D[x]$ have degree $n$, and let $V(p)$ be the set of right roots of $p$ in $D$. For a right root $a$ of $p$, the set 
		\[
		E(p,a)\colonequals \{r\in D\setminus \{0\}\, |\, p(rar^{-1})=0\}\cup \{0\}
		\]
		is a right vector space over the centralizer of $a$ in $D$, which is denoted by $C(a)$. It is known that $V(p)$ intersects at most $n$ conjugacy classes in $D$. Select  $a_1,\dots,a_m\in V(p)$ so that every right root of $p$ is conjugate to exactly one of the elements $a_1,\dots,a_m$. Then the following inequality holds:
		\begin{equation}\label{inequ: equality for roots}
			\sum_{i=1}^m \dim_{C(a_i)} E(p,a_i)\leq n. 
		\end{equation}
	
	\end{subsection}

 \begin{subsection}{Left algebraicity}

	Let $E$ be a division subring of a division ring $D$. 
	An element $a\in D$ is called \textit{left algebraic}  over  $E$  if $a$ is a right root of a nonzero polynomial $p(x)\in E[x]$. It is easy to verify that $a\in D$ is left algebraic over $E$ if and only if the left $E$-vector space $\sum_{n\geq 0}Ea^n$ has a finite dimension. For a left algebraic element $a\in D$ over $E$, the set of all $p(x)\in E[x]$ satisfying $p(a)=0$ forms a left ideal of $E[x]$. Since $E[x]$ is a (left and right) principal ideal ring (PID), this left ideal has a unique monic generator, called the \textit{minimal left polynomial} of $a$ over $E$. The degree of the minimal left polynomial is called the \textit{left degree} of $a$ over $E$. Note that the left degree of $a$ over $E$ coincides with the dimension of the left $E$-vector space $\sum_{n\geq 0}Ea^n$.  We say that a set $S\subset D$ is \textit{left algebraic over $E$ of bounded degree} $n$ if every $a\in S$ is left algebraic over $E$ of left degree $\leq n$. 
	
	Unlike the commutative case, the minimal left polynomial of an element may not necessarily be  irreducible. As an example, the minimal polynomial of $j\in\mathbb{H}$ over $\mathbb{C}=\mathbb{R}+i\mathbb{R}$ is $x^2+1$, which is reducible over $\mathbb{C}$. Nevertheless, we have the following result, which can easily be proved using the product formula \ref{equ: left product formula in one variable}. 
	\begin{proposition}
		Let $a\in D$ be left algebraic over $E$. If the minimal left polynomial of $a$ over $E$ is reducible, then a conjugate of $a$ lies in $E$, i.e. $rar^{-1}\in E$ for some nonzero $r\in D$. 
	\end{proposition} 

\end{subsection}

\begin{subsection}{Right algebraically closed division rings}
	 A division ring $D$ is called \textit{right algebraically closed} (\textit{right AC}) if every nonconstant polynomial $p(x)\in D[x]$ has at least one right root in $D[x]$.  
	The concept of a \textit{left algebraically closed} division ring is defined analogously. It is easy to see that $D$ is right algebraically closed if and only if every irreducible polynomial in the polynomial ring $D[x]$ has degree 1. It follows immediately that a division ring is  right algebraically closed if and only if it is left algebraically closed.

	An example of a right algebraically closed division ring is the division ring of quaternions over $\mathbb{R}$, a result due to Niven and Jacobson. Baer's Theorem states that a noncommutative, centrally finite,  right algebraically closed division ring is isomorphic to the division ring of quaternions over a real-closed field.  Recall that a commutative field $R$ is called \textit{real-closed} if the equation $x^2=-1$ has no solutions in $R$ and the field extension $R(\sqrt{-1})$ is algebraically closed. These results provide a comprehensive characterization of centrally finite, right algebraically closed division rings. By contrast, the case of   centrally infinite, right algebraically closed division rings is not well understood. 
	
	A result of Cohn states that every nonconstant polynomial $p(x)\in D[x]$ over a division ring $D$ has a right root in some extension of $D$ \cite{cohn1975equation}. Using Cohn's result and transfinite induction, one can show that every division ring can be embedded in a right algebraically closed division ring.
	
 \end{subsection}
 
 \end{section}
 
 \begin{section}{Left algebraicity over centralizers}\label{sec: algebraic elements} 
In general, it seems dif and only ificult to make substantial assertions about the structure of the set of left algebraic elements over a division subring. However, the situation becomes more tractable  when the division subring is the centralizer of a given element. As we will see in the next chapter, studying centralizers can help us to understand the class of division rings over which the central Nullstellensatz holds. Therefore, we devote this section to the study of left algebraicity over centralizers. 
\begin{subsection}{Basic results}
  
	Let us fix a division ring $D$ throughout this subsection. Given a subsets $S$ of $D$, we set
	\[
	C(S)\colonequals \{a\in D\,|\, sa=as \text{ for all } s\in S\}.
	\]
	We adopt the convention $C(a_1,\dots,a_n)\colonequals C(\{a_1,\dots.a_n\})$. It is evident that for any subset $S$ of $D$, the centralizer $C(S)$ is a division subring of $D$. To explore algebraicity over centralizers, we begin with the following lemma. 
	\begin{lemma}\label{lem: left-right algebraic}
		Let $b\in D$ be a right root of a nonzero polynomial $p(x)=a_nx^n+\cdots+a_1x+a_0\in D[x]$. Then, $C(a_0,\dots,a_{n})$ is  right algebraic over $C(b)$ with right degree at most $n$. Moreover, the right $C(b)$-vector space generated by $C(a_0,\dots,a_{n})$ is finite-dimensional with dimension at most $n$. 
	\end{lemma} 
	\begin{proof}
		Observe that   	
		\[
		E(p,b)=\{r\in D\,|\,  a_nrb^n+a_{n-1}rb^{n-1}+\cdots+a_1rb+a_0r=0\}.
		\]
		As seen in Subsection \ref{subsec: evaluation}, $E(p,b)$ is a right vector space over $C(b)$ and $\dim_{C(b)} E(p,b)\leq n$ by Inequality \ref{inequ: equality for roots}. Clearly,  $C(a_0,\dots,a_{n})\subset E(p,b)$. In particular, for each $a\in C(a_0,\dots,a_{n})$, the elements $1,a,\dots.,a^{n}$ are right $C(b)$-linearly independent. Hence, $C(a_0,\dots,a_{n})$ is right algebraic over $C(b)$ with right degree at most $n$, which completes the proof. 
	\end{proof}
	
	In what follows, let $S$ be a subgroup of  the multiplicative group $D\setminus \{0\}$. Assume that $b\in D$ is left algebraic over $C(S)$, and let 	
	\[ p(x)=x^n+a_{n-1}x^{n-1}+\cdots+a_1x+a_0\in C(S)[x]\] be the minimal left polynomial of $b$ over $C(S)$. 
	\begin{proposition}\label{prop: left algebraic with subgroup}
		1) $\cap_{r\in S}D[x](x-rbr^{-1})=D[x]p(x)$.\\
		2) Every right root of $p(x)$ in $D$ is conjugate to $b$. \\
		3) The dimension of $E(p,b)$ as a right $C(b)$-vector space is $n$. 
	\end{proposition}
	\begin{proof}
		1) Since $D[x]$ is a principal ideal ring, we have $\cap_{r\in S}D[x](x-rbr^{-1})=D[x]q(x)$ for some monic $q(x)\in D[x]$. Note that $p(x)\in D[x]q(x)$ because $p(rbr^{-1})=0$ for all $r\in S\subset  E(p,b)$. For every $s\in S$, we have
		\[
		s\left( \cap_{r\in S}D[x](x-rbr^{-1}) \right)s^{-1} =D[x]sq(x)s^{-1},
		\]
		implying
		\[
		\cap_{r\in S}D[x](x-srbr^{-1}s^{-1}) =D[x]sq(x)s^{-1}\implies D[x]q(x)=D[x]sq(x)s^{-1}.
		\]
		It follows that $q(x)=sq(x)s^{-1}$ for all $s\in S$, that is, $q(x)\in C(S)[x]$. Since $q(b)=0$, we conclude that $q(x)=p(x)$.\\
		The second statement follows from 1), and 3) follows from 2). 
	\end{proof}
	\begin{remark}
		It follows from the  first part of the proposition that the polynomial $p(x)$ is a Wedderburn polynomial, a concept defined and studied  by Lam and Leroy \cite{lam2004wedderburn}. 
	\end{remark}
\end{subsection}
\begin{subsection}{Algebraicity over $C(a)$}
	As before, let $D$ be a fixed division ring throughout this subsection. 
	\begin{proposition}\label{prop: leftalg-rightalg}
		For any $a,b\in D$, $b$ is left algebraic over $C(a)$ if and only if $a$ is right algebraic over $C(b)$, in which case the left degree of $b$ over $C(a)$ coincides with  the right degree of $a$ over $C(b)$. 
	\end{proposition}
	\begin{proof}
		Let \[ p(x)=x^n+a_{n-1}x^{n-1}+\cdots+a_1x+a_0\] be the minimal left polynomial of $b$ over $C(a)$. By Lemma \ref{lem: left-right algebraic}, every element of $C(a_0,\dots,a_{n-1})$ is right algebraic over $C(b)$ with right degree at most  $n$. Since $a\in C(a_0,\dots,a_{n-1})$, it follows that $a$ is  right algebraic over $C(b)$ with right degree at most $n$. By symmetry, the left degree of $b$ over $C(a)$ equals the right degree of $a$ over $C(b)$.
	\end{proof}
	As a consequence, we obtain the following description of the set of left algebraic elements over the centralizer of an element. 
	\begin{corollary}\label{cor: closure of C(a)}
		The set of all left algebraic elements of $D$ over $C(a)$ is 
		$\cup C(a_0,\dots,a_{n-1})$,
		where the union ranges over all $(a_0,\dots,a_{n-1})\in D^{n}$ satisfying
		\[
		a^n+a^{n-1}a_{n-1}+\cdots+aa_1+a_0=0.
		\]
	\end{corollary}
	\begin{proof}
		Let  $(a_0,\dots,a_{n-1})\in D^{n}$ satisfy the equation
		\[
			a^n+a^{n-1}a_{n-1}+\cdots+aa_1+a_0=0.
		\]
		It follows from Lemma \ref{lem: left-right algebraic} that every element of $C(a_0,\dots,a_{n-1})$ is left algebraic over $C(a)$. This proves that the union is left algebraic over $C(a)$. Now, let $b\in D$ be 
		left algebraic over $C(a)$. By Proposition \ref{prop: leftalg-rightalg}, $a$ is right algebraic over $C(b)$. Let   \[ p(x)=x^n+x^{n-1}a_{n-1}+\cdots+xa_1+a_0\]
		be the minimal right polynomial of $a$ over $C(b)$. We have
		\[
			0=p_r(a)=a^n+a^{n-1}a_{n-1}+\cdots+aa_1+a_0.
		\] 
		Since $b\in C(a_0,\dots,a_{n-1})$, the result follows. 
		
	\end{proof}

\end{subsection}
\begin{subsection}{A rational criterion for algebraicity over $C(a)$}
Let $D$ be a division ring, and let $F(x_0,x_1,\dots)$ denote the set of all rational expressions in $x_0,x_1,\dots$ with coefficients from the prime field $F$ of $D$. For a comprehensive discussion of rational expressions, the reader is referred to Cohn's book \cite[Subsecion 7.2]{cohn1995skew}. Inductively, we define a sequence of rational expressions $L_1,L_2,\dots$ in $F(x_0,x_1,\dots)$  as follows: $L_1(x_0,x_1)\colonequals x_1$, and for $n\geq 1$, 
\[
	L_{n+1}(x_0,x_1,\dots,x_{n+1}) \colonequals L_n(x_0,[x_0,x_1x_{n+1}^{-1}],\dots,[x_0,x_nx_{n+1}^{-1}]),
\] 
where $[a,b]\colonequals ab-ba$.
\begin{lemma}\label{lem: linearly independence}
 For $a,b_1,\dots,b_n\in D$, the elements $b_1,\dots,b_n\in D$ are left linearly independent over $C(a)$ if and only if $L_n(a,b_1,\dots,b_n)$ is defined and nonzero.    
\end{lemma}
\begin{proof}
	The key observation is that $b_1,\dots,b_n\in D$ are left linearly independent over $C(a)$ if and only if $b_n\neq 0$ and $[a,b_1b_{n}^{-1}],\dots,[a,b_{n-1}b_{n}^{-1}]$ are left linearly independent over $C(a)$. A straightforward  induction on $n$ then proves the result.
\end{proof}

\begin{remark}
	It is easy to see that $L_n(a,b_1,\dots,b_n)$ is zero if and only if  $b_1,\dots,b_n$ are left linearly dependent over $C(a)$, while every proper subset of  $\{b_1,\dots,b_n\}$ remains left linearly dependent over $C(a)$. 
\end{remark}

By defining $F_n(x,y)\colonequals L_n(x,1,y,\dots,y^n)$, we can derive a rational criterion for left algebraicity over the centralizer of an element:  
\begin{proposition}\label{prop: rational identity for algebraicity}
	For any $a,b\in D$, the element $b$ is left algebraic of left degree $n$ over $C(a)$ if and only if $F_{n}(a,b)=0$.  
\end{proposition}
\begin{proof}
	Note that $b$ is left algebraic of left degree $n$ over $C(a)$ if and only if $1,b,\dots,b^{n-1}$ are left linearly independent over $C(a)$, while
	$1,b,\dots,b^{n}$ are left linearly dependent over $C(a)$. An application of Lemma \ref{lem: linearly independence} gives the result. 
\end{proof}
As an application of of this criterion, we present the following result. 
\begin{proposition}
	Let $D$ be  a division ring with infinite center. Suppose there exists a natural number $n$ such that every $b\in D$ is left algebraic of left degree $\leq n$ over $C(a)$ for all $a\in D$. Then, $D$ is centrally finite. 
\end{proposition}
\begin{proof}
	By Proposition \ref{prop: rational identity for algebraicity}, the rational expression $F_n(x,y)$ is a rational identity for $D$. Since the center of $D$ is infinite and $D$ satisfies a (nontrivial) rational identity, it follows from foundational results in the theory of generalized rational identities (see \cite{rowen1980polynomial}) that $D$ is centrally finite.   
\end{proof}
Clearly, the converse of the proposition holds true. 
\begin{remark}
	It is worth noting that Bell et al. have proved that if a division ring  is of bounded left degree over a (commutative) subfield, then it is centrally finite \cite{bell2013shirshov}.
\end{remark}

\end{subsection}

\end{section} 


\begin{section}{The central Nullstellensatz}\label{sec: null}
This section contains the main result of this article. More precisely, we establish necessary and sufficient conditions under which a division ring satisfies the central Nullstellensatz. 
\begin{subsection}{Centrally algebraically closed division rings} 
 Let $D$ be a division ring, and define $D^n_c$ as the set of all $(a_1,\dots,a_n) \in D^n$ such that $a_ia_j=a_ja_i$ for all  $i,j$.  
 
 \begin{definition}
 	A division ring $D$ is said to be centrally algebraically closed (\textit{centrally AC}) if, for every $(a_1,\dots,a_n) \in D^n_c$, the centralizer $C(a_1,\dots,a_n)$ is right algebraically closed. 
 \end{definition}
We record the following two basic facts regarding centrally AC division rings. 
\begin{proposition}
	If $D$  is a centrally AC division ring, then for any point $(a_1,\dots,a_n) \in D^n_c$, the centralizer $C(a_1,\dots,a_n)$ is itself centrally AC.
\end{proposition}
\begin{proof}
	Trivial! 
\end{proof}
\begin{proposition}\label{prop: quaternions centrally AC}
	The division ring of quaternions over a real-closed field is centrally AC.
\end{proposition}
\begin{proof}
	Let $\mathbb{H}$ be the division ring of quaternions over a real-closed field  $R$. It is easy to see that  for any $(a_1,\dots,a_n) \in \mathbb{H}^n_c$, the centralizer $C(a_1,\dots,a_n)$ is either $\mathbb{H}$  or the commutative field $R(\sqrt{-1})$, both of which are right algebraically closed.
\end{proof}
\end{subsection}

\begin{subsection}{The central Nullstellensatz over division rings}
	We now come to the main result of the paper. 
	\begin{theorem}\label{thm: nullstelensatz for n}
	Let $D$ be a division ring and $n\geq 1$. Then the following conditions are equivalent: \\
	(1) The set of all maximal left ideals of $R_n\colonequals D[x_1,\dots,x_n]$ coincides with the set of all left ideals of the form 
	$$R_n(x_1-a_1)+\cdots+R_n(x_n-a_n),$$ where $(a_1,\dots,a_n)\in D_c^n$.\\
	(2) Any simple left $D[x_1,\dots,x_n]$-module is one-dimensional as a left $D$-vector space.\\
	(3) For any $(a_1,\dots,a_{n-1})\in D^n_c$, the centralizer $C(a_1,\dots,a_{n-1})$ is right algebraically closed. In the case $n=1$, the condition reduces to $D$ being right algebraically closed. \\
	(4) The set of all maximal right ideals of $R_n$ coincides with the set of all right ideals of the form 
	$$(x_1-a_1)R_n+\cdots+(x_n-a_n)R_n,$$ where $(a_1,\dots,a_n)\in D_c^n$.\\ 	
	\end{theorem}	
	\begin{proof}
	It is easy to prove the implication (1)$\implies$(2), and the implication (2)$\implies$(1) is proved in \cite[Lemma 2.2]{chapman2025amitsur}. It is also easy to show that the implication (2)$\implies$(3) holds for $n=1$. To prove the implication 
	(2)$\implies$(3) for $n>1$, we use proof by contradiction. Let $(a_1,\dots,a_{n-1})\in D^n_c$ be a point such that $C(a_1,\dots,a_{n-1})$ is not right algebraically closed. It follows that there exists a nonconstant polynomial $p(x)\in C(a_1,\dots,a_{n-1})[x] $ that is irreducible over $C(a_1,\dots,a_{n-1})$. 
	Consider the left $D$-vector space $M=D[x]/D[x]p$. It is easy to see that $M$ can be turned into a left $D[x_1,\dots,x_n]$-module using
	\[\begin{cases}
		x_i(q+D[x]p)=qa_i+D[x]p & \text{ for } i<n,\\
		x_n(q+D[x]p)=xq+D[x]p. 
	\end{cases} \] 
	 Since $M$ is finite-dimensional as a left $D$-vector space, $M$ contains a simple $D[x_1,\dots,x_n]$-submodule $W$. By (2), $W=D(q_0+D[x]p)$ for some monic polynomial $q_0\in D[x]$ of degree less than $\deg(p)$. Note that $q_0\in C(a_1,\dots,a_{n-1})[x]$ because $q_0a_i\in Dq_0$ for $i<n$ and $q_0$ is monic. Since $x_n(q_0+D[x]p)=b(q_0+D[x]p)$ for some $b\in D$, we conclude that $(x-b)q_0=p$, which contradicts the assumption that $p$ is irreducible over $C(a_1,\dots,a_{n-1})$.

	We prove the implication (3)$\implies$(2) by induction on $n$. The base case $n=1$ is immediate since $D[x_1]$ is a principal ideal domain. Now suppose $n>1$, and let $M$ be a simple left $R_n$-module. By the Amitsur--Small Theorem \ref{thm: amitsur small}, $M$ is finite-dimensional as a $D$-vector space. By induction, $M$ contains be a nonzero $R_{n-1}$-submodule $W=Dv_0$ for some $v_0\in M$. Consequently, for each $i<n$, there exists $a_i\in D$ such that $x_iv_0=a_iv_0$. It is easy to straightforward that  $(a_1,\dots,a_{n-1})\in D^n_c$. Since $M$ is finite-dimensional as a $D$-vector space and $D[x]$ is a (left and right) principal ideal domain, the set
	\[
		\{P(x)\in D[x]\, |\, P(x_n)v_0=0\}
	\]
	is a principal left ideal generated by a monic polynomial $p(x)\in D[x]$. For each $i<n$,  
	\[
		0=x_ip(x_n)v_0=p(x_n)x_iv_0=p(x_n)a_iv_0.
	\] By the uniqueness of $p$, it follows that $pa_i=a_ip$, that is, $p\in C(a_1,\dots,a_{n-1})[x]$. Since $C(a_1,\dots,a_{n-1})$ is right (and left) algebraically closed, $C(a_1,\dots,a_{n-1})$ contains a  left root $a$ of $p$. Hence, $p=(x-a)q$ for some polynomial $q\in C(a_1,\dots,a_{n-1})[x]$. The reader can verify that $Dq(x_n)v_0$ is an $R_n$-submodule of $M$. Since $M$ is simple, it follows that $M=Dq(x_n)v_0$. In particular, $M$ is one-dimensional as a $D$-vector space, completing the proof of (3)$\implies$(2). The equivalence of (4) and (3) follows from symmetry. 
	\end{proof}
	Following Chapman and Paran, we call a division ring a \textit{Nullstellensatz ring} if the central Nullstellensatz holds for the division ring \cite[Definition 4.1]{chapman2025amitsur}. A direct consequence of Theorem \ref{thm: nullstelensatz for n} is the following characterization. 
	\begin{theorem}\label{thm: nullstelensatz}
		A division ring is a Nullstellensatz ring if and only if it is centrally algebraically closed. 	
	\end{theorem}

	In the light of Proposition \ref{prop: quaternions centrally AC}, our discussion provides a proof of Theorem \ref{thm: Weak Nullstellensatz over quaternions}, which is fundamentally different from the proof given by Alon and Paran in \cite{AlonNullestellensatz}. It is noteworthy  that  the theorem implies that the notion of a Nullstellensatz ring is left-right symmetric, as is the concept of a centrally AC division ring. 
		
	In \cite[Definition 3.1]{chapman2025amitsur}, the authors defined the notion of an Amitsur--Samll ring as a division ring $D$ that satisfies the following property: For every maximal left ideal $I$ of $D[x_1,\dots,x_n]$, the left ideal of $I\cap D[x_1,\dots,x_m]$ is maximal for all $m<n$. They proved that a right AC division ring is an Amitsur--Samll ring  if and only if it is a Nullstellensatz ring \cite[Theorem 4.2]{chapman2025amitsur}. We end this part with a more complete version of their result:
	\begin{proposition}
			For a right AC division ring, the following statements are equivalent:\\
			1) $D$ is a Nullstellensatz ring.\\
			2) $D$ is an Amitsur--Samll ring.\\
			3) $D$ is centrally AC. 
	\end{proposition} 
		
\end{subsection}
\begin{subsection}{Hilbert's Nullstellensatz over centrally algebraically closed division rings}
	Let us begin with some preliminaries. For further details, the reader is referred to \cite{aryapoor2024explicit} and the references therein. Let $D[x_1,\dots,x_n]$ denote the ring of polynomials in $n$ central indeterminates over a division ring  $D$. Every polynomial in $D[x_1,\dots,x_n]$ can be evaluated on $D^n_c$ as a left polynomial via substitution. Give a subset $S$ of $D[x_1,\dots,x_n]$, the c-zero locus of $S$ is defined to be the set of all points $(a_1,\dots,a_n)\in D^n_c$ that are annihilated by every polynomial in $S$. The following version of Hilbert's Nullstellensatz is valid for centrally AC division rings.  
	\begin{theorem}[Hilbert's Nullstellensatz]\label{thm: Explicit Nullstellensatz}
	Let $D$ be a centrally algebraically closed division ring. A polynomial $p\in D[x_1,\dots, x_n]$ vanishes on the c-zero locus of a left ideal $I$ of $D[x_1,\dots, x_n]$  if and only if for every $a\in D$, there exists a natural number $N\geq 1$ such that 
	\[(ap)^N\in I+I(ap)+I(ap)^2+\cdots+I(ap)^N.\]
	\end{theorem}
	\begin{proof}
		The same argument as in \cite[Theorem 1.5]{aryapoor2024explicit} for the case over $\mathbb{H}$ applies \textit{mutatis mutandis} in this setting. Observe that the proof in \cite{aryapoor2024explicit} relies on the central Nullstellensatz and the Rabinowitsch trick, both of which remain valid over any centrally AC division ring.   
	\end{proof}
 
\end{subsection}
\end{section} 


\begin{section}{Further results and open problems}
In the final section of the paper, we present more results regarding centrally AC division rings, and moreover, present some natural questions. 
\begin{subsection}{Existentially closed division rings}
 
 We begin by recalling the definition of an existentially closed division ring (for more details, see Cohn's book \cite{cohn1995skew}). Let $D$ be a division ring with center $k$, and  let $D_k\langle x_1,\dots,x_n\rangle$ denote the ring of polynomials in $n$ non-commuting indeterminates subject to the relations $ax=xa$, where $a\in k$. The division ring $D$  is said to be  \textit{existentially closed} (\textit{EC}) if  any system of equations $f_1=\cdots=f_m=0$, where $f_1,\dots,f_m\in D_k\langle x_1,\dots,x_n\rangle$, that admits a solution in some extension of $D$, already has  a solution in $D$.   
\begin{proposition}
	Every EC division ring is necessarily centrally AC. 
\end{proposition}
\begin{proof}
	We first note that every EC division ring is right AC since any left polynomial over a division ring has a right solution in some extension of the division ring \cite{cohn1975equation}. The result now follows from the fact that the centralizer of an element in an EC division ring is also EC \cite[Thoerem 6.5.8]{cohn1995skew}. 
\end{proof}
To the best of the author's knowledge, no explicit example of an EC division ring is known. It is, however, known that every division ring $D$ can be embedded into an EC division ring that shares the same center \cite[Theorem 6.5.3]{cohn1995skew}, from which the following result follows immediately.
\begin{theorem}
	Every division ring admits an extension (with the same center) that is centrally AC. 
\end{theorem}

\end{subsection}	
\begin{subsection}{Questions}

	We conclude the paper with some open questions concerning centrally AC division rings:\\
	\textbf{Question 1}: Does there exist a right algebraically closed division ring that is not centrally algebraically closed? \\
	\textbf{Question 2}: Is every characteristically algebraically closed division ring also centrally AC? Recall that a division ring $D$ is called \textit{characteristically algebraically closed}  if for every square matrix $A$ over $D$, there exists an element $a\in D$ such that $A-aI$ is not invertible.   In \cite{wood1985quaternionic}, Wood proved that the ring of quaternions over $\mathbb{R}$ is characteristically algebraically closed. \\
	\textbf{Question 3}: Is every polynomially algebraically closed division ring necessarily centrally AC? A division ring $D$ with center $k$ is called \textit{polynomially algebraically closed} (\textit{polynomially AC}) if every nonconstant polynomial in $ D_k\langle x\rangle$ has a root in $D$. The first example of a polynomially AC division ring was constructed by Makar-Limanov \cite{makar1985algebraically}.
	 
\end{subsection}

 
\end{section} 


\bibliographystyle{plain}
\bibliography{ndbiblan}

 \end{document}